\documentstyle[11pt,dgspp]{article}
\newtheorem{theorem}{Theorem}[section]
\newtheorem{definition}[theorem]{Definition}
\newtheorem{proposition}[theorem]{Proposition}
\newtheorem{lemma}[theorem]{Lemma}
\newtheorem{corollary}[theorem]{Corollary}
\newremark{example}[theorem]{Example}
\newremark{remark}[theorem]{Remark}

\newcommand{\bds}{\bigoplus} %big direct sum

\newcommand{\bul}{\vcenter{\hbox{$\scriptscriptstyle\,\bullet$}}\,}

\newcommand{\da}{\dot\alpha}

\newcommand{\ddv}{\ddot{v}}

\newcommand{\del}[1]{\nabla_{#1}}
\newcommand{\dg}{\dot\gamma}

\newcommand{\dps}{\displaystyle}
\newcommand{\ds}{\oplus} %direct sum

\newcommand{\dv}{\dot{v}}

\newcommand{\dz}{\dot{z}}

\newcommand{\g}{\gamma}

\newcommand{\half}{\mbox{$\txs\frac{1}{2}$}}

\newcommand{\io}{\iota}
\newcommand{\Iaut}{I^{a\kern-.05em u\kern-.05em t}}
\newcommand{\Ispl}{I^{s\kern-.055em p\kern-.025em l}}
\newcommand{\jac}{R_{\,\dg}}
\newcommand{\k}{\frak{k}}
\newcommand{\lsp}{[\kern-0.15em[} %left spanning bracket

\newcommand{\n}{\mbox{$\frak{n}$}}

\newcommand{\ph}{$p\kern-.1em H\!$}
\newcommand{\quar}{\mbox{$\txs\frac{1}{4}$}}
\newcommand{\rsp}{]\kern-0.15em]} %right spanning bracket

\newcommand{\sixth}{\mbox{$\txs\frac{1}{6}$}}

\newcommand{\surj}{\rightarrow\kern-.82em\rightarrow}

\newcommand{\tquar}{\mbox{$\txs\frac{3}{4}$}}

\newcommand{\txs}{\textstyle}
\newcommand{\ve}{\mbox{$\varepsilon$}}

\newcommand{\w}{\frak{w}}

\newcommand{\z}{\frak{z}}
\newcommand{\E}{\frak{E}}

\newcommand{\II}{\rscr{I{\kern-.55em}I}}

\newcommand{\R}{\Bbb R}
\newcommand{\U}{\frak{U}}
\newcommand{\V}{\frak{V}}
\newcommand{\Z}{\frak{Z}}
\newcommand{\bZ}{\Bbb Z}

\renewcommand{\l}{\lambda}

\renewcommand{\v}{\frak v}

\makeatletter
\newcommand{\Ad}[1]{\mathop{\operator@font Ad}\nolimits_{#1}}
\newcommand{\ad}[1]{\mathop{\operator@font ad}\nolimits_{#1}}
\newcommand{\add}[2]{\mathop{\operator@font
   ad}\nolimits^{\dagger}_{#1}{\!#2}}
\newcommand{\Aut}{\mathop{\operator@font Aut}\nolimits}
\newcommand{\diag}{\mathop{\operator@font diag}\nolimits}
\newcommand{\End}{\mathop{\operator@font End}\nolimits}
\newcommand{\im}{\mathop{\operator@font im}\nolimits}
\newcommand{\mcp}{\mathop{\operator@font mult}\nolimits_{cp}}
\newcommand{\mev}{\mathop{\operator@font mult}\nolimits_{ev}}
\newcommand{\Ric}{\mathop{\operator@font Ric}\nolimits}
\newcommand{\Dspecp}{\mathop{\operator@font {\escr D}spec}\nolimits_\wp}
\newcommand{\specl}{\mathop{\operator@font spec}\nolimits_\ell}
\newcommand{\specp}{\mathop{\operator@font spec}\nolimits_\wp}
\newcommand{\tr}{\mathop{\operator@font tr}\nolimits}
\makeatother

\preprint{JP1}
\title{Conjugate Loci of Pseudoriemannian 2-step Nilpotent Lie
Groups with Nondegenerate Center}
\author{Changrim Jang\thanks{On sabbatical leave; permanent address: 
        Department of Mathematics,
        College of Natural Sciences,
        University of Ulsan,
        Ulsan 680-749,
        Republic of Korea.\protect\newline Email:
        crjang@mail.ulsan.ac.kr}
\quad {\small and} \quad Phillip E. Parker}
\address{Mathematics Department\\
        Wichita State University\\
        Wichita KS 67260-0033\\
        USA\\
        \hspace*{-1em}crjang@math.wichita.edu \qquad\quad 
        phil@math.wichita.edu}
\date{23 March 2005}
\abstract{
We determine the complete conjugate locus along all geodesics parallel or
perpendicular to the center (Theorem \ref{th1}).  When the center is
1-dimensional we obtain formulas in all cases (Theorem \ref{th2}), and
when a certain operator is also diagonalizable these formulas become
completely explicit (Corollary \ref{cor}).  These yield some new
information about the smoothness of the pseudoriemannian conjugate locus.
We also obtain the multiplicities of all conjugate points.
}
\msc{53C50}{22E25, 53B30, 53C30}

\begin{document}

\maketitle

% uncomment to print 1-sided title page and rest 2-sided
%\setcounter{page}{0}\thispagestyle{empty}\strut\vfill\eject

\section{Introduction}
\label{intro}

Let $N$ denote a nilpotent Lie group with a left-invariant metric tensor. 
The study of the conjugate locus in these spaces began with O'Sullivan's 
proof in 1974 that when $N$ is not abelian, there exist conjugate points 
\cite{OS}.  In 1981, Naitoh and Sakane showed \cite{NS} that some first
conjugate points must lie in the center $Z$ of $N$.  If $N$ is 2-step
nilpotent, they showed that there are no closed geodesics and the cut
locus of the identity element $1\in N$ is nonempty.  Finally, when $N$ is
the 3-dimensional Heisenberg group $H_3$, the cut locus of 1 is precisely
the set of first conjugate points to 1 and lies in the center of $H_3$.
Since then, all results that have been published were for 2-step 
nilpotent Lie groups. This class includes the groups of $H$-type; see the 
introduction of \cite{JP4} for a brief history of the latter.

Eberlein's 1994 paper \cite{E} sparked a continuing increased interest in
2-step nilpotent Lie groups.  As he noted, they provide examples in which
explicit calculations are frequently feasible and which differ minimally
from flat (pseudo-) Euclidean spaces.  He obtained only one result on
conjugate points, determining conditions under which geodesics tangent to
the center contain conjugate points.  So far, only Riemannian
(positive-definite) metric tensors had been considered.  Cordero and
Parker \cite{CP4} extended much of this work to pseudoriemannian
(indefinite) metric tensors, including the result on conjugate points; 
{\em cf.}\ Proposition \ref{4.18} {\em infra}.

Meanwhile, in 1997, Walschap \cite{W} showed that for a certain class of
such Riemannian groups (nonsingular with 1-dimensional center), the cut
locus and the conjugate locus coincide, and he made an explicit
determination of all first conjugate points in them.  Jang and Park
\cite{JPk} later gave explicit formulas for all conjugate points on
geodesics either parallel or perpendicular to the center in any Riemannian
2-step nilpotent Lie group, and on all geodesics when the center is
1-dimensional. The present paper extends these results to pseudoriemannian
geometries, while providing shorter, simpler, and more conceptual proofs. 
We also obtain some new information about the smoothness of the conjugate 
locus which are valid even when there are degenerate \cite{MPT} conjugate 
points, which is precisely when the well-known results of Warner 
\cite{Wa} do not apply.  Indeed, Warner's property (R2) is equivalent to
conjugate points being nondegenerate, as is easily seen from his
interpretation of it on p.\,577f (and proof on p.\,603) of \cite{Wa}.

\medskip
By an {\em inner product\/} on a vector space $V$ we shall mean a
nondegenerate, symmetric bilinear form on $V$, generally denoted by
$\langle\,,\rangle$.  Our convention is that $v\in V$ is timelike if
$\langle v,v\rangle > 0$, null if $\langle v,v\rangle = 0$, and spacelike
if $\langle v,v\rangle < 0$.

Throughout, $N$ will denote a connected, simply connected, 2-step
nilpotent Lie group with Lie algebra \n\ having center $\z$.  This means 
the Lie group exponential map $\exp :\n\to N$ is a global diffeomorphism, 
and we may work in what are usually called ``exponential coordinates," 
making calculations in $\n$ and taking the result to be in $N$.  We shall
use $\langle\,,\rangle$ to denote either an inner product on \n\ or the
induced left-invariant pseudoriemannian (indefinite) metric tensor on $N$.
Note that we shall {\em not\/} be using the geodesic (pseudoriemannian) 
exponential map in this article.

The only general result on conjugate points is Proposition 4.18 from
\cite{CP4}; {\em cf.}\ Proposition 3.10 in Eberlein \cite{E}.
\begin{proposition}
Let\label{4.18} $N$ be a simply connected, 2-step nilpotent Lie group 
with left-invariant metric tensor $\langle\,,\rangle$, and let $\g$ be a
geodesic with $\dg(0) = a\in \z$.  If\/ $\add{\bul}{a} = 0$, then there
are no conjugate points along $\g$.
\end{proposition}

In this paper we make the additional assumption that the center $\z$ of
$\n$ is nondegenerate.  (In the notation of \cite{CP4}, this means $\U =
\V = \{0\}$, $\Z = \z$, $\E = \v = \z^\perp$ and $\n = \z\ds\v$.)

Theorems \ref{th1} and \ref{th2} and Corollary \ref{cor} generalize and
extend the results of \cite{JPk} to pseudoriemannian spaces.  See also
\cite{PT} for results on the topology of the set of conjugate points
along a single geodesic.

Corollaries \ref{avz} and \ref{smz} and Proposition \ref{scl} give some
new information on the smoothness of the conjugate locus.

Finally, we also obtain a new formulation of the Jacobi equation for these 
geometries in Proposition \ref{yj}, of some interest in its own right.

\medskip
We denote the adjoint with respect to $\langle\,,\rangle$ of the adjoint
representation of the Lie algebra \n\ on itself by $\add{}{}\!$.  In the
case of a nondegenerate center, the involution $\iota$ of \cite{CP4} is
merely given by $\io (z_\alpha )=\ve_\alpha\, z_\alpha$ and $\io (e_a )=
\bar{\ve}_a \, e_a$ where $\langle z_\alpha ,z_\alpha \rangle =
\ve_\alpha$ and $\langle e_a,e_a \rangle = \bar{\ve}_a$ on an orthonormal
basis of \n.  Then the operator $j:\z\rightarrow\End\left(\v\right) $ is
given by $ j(z)x = \io\add{x}{\io z}$.  We refer to \cite{CP4} and 
\cite{O} for other notations and results.

For convenience, we shall use the notation $J_z = \add{\bul}{z}$ for any
$z\in\z$.  (The involution $\io$ does not play much of a role when the
center is nondegenerate, and in fact may be omitted.)

We begin by specializing Theorems 3.1 and 3.6 of \cite{CP4} to
the case of a nondegenerate center.
\begin{proposition}
For\label{conn} all $z,z'\in\z$ and $e,e'\in\v$ we have
\begin{eqnarray*}
\del{z}z' & = & 0\, ,\\
\del{z}e = \del{e}z &=& -\half J_z e\, ,\\
\del{e}e' &=& \half [e,e']\,.
\end{eqnarray*}
\end{proposition}
\begin{proposition}
For\label{curv} all $z,z',z''\in\z$ and $e,e',e''\in\v$ we have
\begin{eqnarray*}
R(z,z')z'' &=& 0\, ,\\
R(z,z')e   &=& \quar (J_z J_{z'} e - J_{z'} J_z e)\, ,\\
R(z,e)z'   &=& \quar J_z J_{z'} e\, ,\\
R(z,e)e'   &=& \quar [e, J_z e']\, ,\\
R(e,e')z   &=& -\quar\big( [e, J_z e'] +[ J_z e,e']\big)\, ,\\
R(e,e')e'' &=& \quar\big( J_{[e,e'']}e' - J_{[e',e'']}e\big)
               +\half J_{[e,e']}e''\,.
\end{eqnarray*}
\end{proposition}

To study conjugate points, we shall use the Jacobi operator.
\begin{definition}
Along\label{jo} the geodesic $\g$, the\/ {\em Jacobi operator} is given by
$$ \jac \bul = R(\bul,\dg)\dg\,. $$
\end{definition}
In physics, this operator measures the relative acceleration produced by
tidal forces along $\g$ \cite[p.\,219]{O}.  For the reader's convenience,
we recall that a {\em Jacobi field\/} along $\g$ is a vector field along
$\g$ which is a solution of the {\em Jacobi equation}
$$ \nabla_{\dg}^2 Y(t)+\jac Y(t) = 0 $$
along $\g$. The point $\g(t_0)$ is {\em conjugate\/} to the point $\g(0)$ 
if and only if there exists a nontrivial Jacobi field $Y$ along $\g$ such 
that $Y(0) = Y(t_0) = 0$.

Next, we specialize Proposition 4.8 of \cite{CP4} to our present setting.
As there, $L_n$ denotes left translation in $N$ by $n\in N$.
\begin{proposition}
Let\label{e3.2} $N$ be simply connected, $\g$ a geodesic with $\g(0)= 
1\in N$, and $\dot{\g}(0) = z_0 + x_0 \in \n$. Then
$$\dot{\g}(t) = L_{\gamma(t)*}\left( z_0 + e^{tJ}x_0 \right) $$
where $J = J_{z_0}$.
\end{proposition}
As in \cite{JPk}, we shall identify all tangent spaces with $\n = T_1 N$.
Thus we regard
$$ \dg(t) = z_0 + e^{tJ}x_0 $$
as the geodesic equation. Using Proposition \ref{curv}, we compute
\begin{lemma}
The\label{jof} Jacobi operator along the geodesic $\g$ in $N$ with 
$\g(0) = 1$ and $\dg(0) = z_0 + x_0$ is given by
\begin{eqnarray*}
\jac(z+x) &=& \tquar J_{[x,x']}x' + \half J_z Jx' - \quar
   JJ_z x' - \quar J^2 x \\[.5ex]
&& {} - \half[x, Jx'] + \quar[x',Jx] + \quar[x',J_z x']
\end{eqnarray*}
for all $z\in\z$ and $x\in\v$, where $x' = e^{tJ}x_0$ and $J=J_{z_0}$.\eop
\end{lemma}

\noindent
{\sc Acknowledgements:} CJ wishes to express his profound gratitude to 
Prof.\ K. Park, who may not have been present {\em in corpore\/} but was 
most definitely present {\em in spiritu\/} during the initial formulation
of some proofs.  CJ also thanks the Department at Wichita for its
hospitality.  PEP thanks CJ for his inexhaustible energy and his patience.

\section{Main results}
\label{ndc}

For later convenience, we separate out the principal part of the necessary
computations for Jacobi fields. This yields a new formulation (due to CJ) 
of the Jacobi equation for these geometries.
\begin{proposition}
Let\label{yj} $\g$ be a geodesic with $\g(0) =1$ and $\dg(0)
= z_0 + x_0 \in \z\ds\v=\n$.  A vector field $Y(t) = z(t) + e^{tJ}v(t)$
along $\g$, where $ z(t) \in \z $ and $v(t) \in \v$ for each $t$, is a
Jacobi field if and only if
\begin{eqnarray*}
\dot{z}(t) -[e^{tJ}v(t) , x'(t)] &=& \zeta\, , \\
e^{tJ}\ddot{v}(t)+e^{tJ}J\dot{v}(t)-J_{\zeta}x'(t) &=& 0\, ,
\end{eqnarray*}
where $x'(t) = e^{tJ}x_0$ with $J=J_{z_0}$ and $\zeta\in\z$ is a constant.
\end{proposition}
\begin{proof}
Using Proposition \ref{conn} and Lemma \ref{jof}, we compute
$\nabla_{\dg}Y, \nabla_{\dg}^2 Y$ and $\jac Y$ for such a $Y$.
\begin{eqnarray*}
\nabla_{\dot{\gamma}} Y(t) 
&=& \dot{z}(t) -\half  J_{z(t)}x'(t) +J e^{tJ}v(t)
   +e^{tJ}\dot{v}(t) \\
&& {}-\half Je^{tJ}v(t) + \half  [x'(t),e^{tJ}v(t)]\\ 
&=&\dot{z}(t)-\half J_{z(t)}x'(t) +\half Je^{tJ}v(t)
   +e^{tJ}\dot{v}(t) + \half  [x'(t), e^{tJ}v(t)] ,
\end{eqnarray*}
\begin{eqnarray*}
\nabla_{\dot{\gamma}}^2 Y(t) 
&=& \ddot{z}(t) -\half J_{\dz(t)}x'(t)-\half  
   J_{\dz(t)}x'(t) - \half J_{z(t)}Jx'(t)\\
&&{} +\quar JJ_{\dz(t)}x'(t)-\quar[x'(t),J_{z(t)}x'(t)]
   +\half J^2 e^{tJ}v(t) \\
&&{} +\half Je^{tJ}\dot{v}(t) - \quar J^2e^{tJ}v(t) + \quar
   [x'(t),Je^{tJ}v(t)] + Je^{tJ}\dot{v}(t) \\
&&{} +e^{tJ}\ddv(t) - \half Je^{tJ} \dot{v}(t) + \half [x'(t),
   e^{tJ}\dot{v}(t)] + \half [J x'(t),e^{tJ}v(t)] \\
&&{} + \half[x'(t), Je^{tJ}v(t) + e^{tJ}\dv(t)] - \quar 
   \add{x'(t)}{[x'(t), e^{tJ}v(t)]} \\
&=& \ddot{z}(t) - J_{\dz(t)}x'(t) - \half J_{z(t)}Jx'(t)\\
&&{} +\quar JJ_{\dz(t)}x'(t) - \quar
   [x'(t),J_{z(t)}x'(t)] + \quar J^2 e^{tJ}v(t) +
   Je^{tJ}\dot{v}(t) \\
&&{} + \tquar [x'(t),Je^{tJ}v(t)] + e^{tJ}\ddot{v}(t) 
   + [x'(t), e^{tJ}\dot{v}(t)] \\
&&{} + \half [J x'(t),e^{tJ}v(t)] - \quar\add{x'(t)}{[x'(t), 
   e^{tJ}v(t)]}\, ,
\end{eqnarray*}
\begin{eqnarray*}
\jac Y(t) &=& \jac z(t) + \jac e^{tJ}v(t)\\
&=& \half J_{z(t)}Jx'(t) - \quar JJ_{z(t)}x'(t) +
   \quar[x'(t),J_{z(t)}x'(t)]\\
&& {}+\tquar\add{x'(t)}{[e^{tJ}v(t),x'(t)]} - \quar J^2 
   e^{tJ}v(t) - \half [e^{tJ}v(t), Jx'(t)]\\
&& {}+\quar[x'(t), Je^{tJ}v(t)]
\end{eqnarray*}
Thus $Y$ is a Jacobi field if and only if
\begin{eqnarray*}
\nabla_{\dg}^2 Y(t)+\jac Y(t)
&=& \ddot{z}(t) - J_{\dz(t)}x'(t) + Je^{tJ}\dot{v}(t) +
   [x'(t),Je^{tJ}v(t)]\\
&& {}+e^{tJ}\ddot{v}(t)+[x'(t),e^{tJ}\dot{v}(t)]+[Jx'(t),e^{tJ}v(t)]\\
&& {}+\add{x'(t)}{[e^{tJ}v(t),x'(t)]}\\
&=& \frac{d}{dt}\Bigl(\dot{z}(t) - [e^{tJ}v(t),x'(t)]\Bigr)  
   + e^{tJ}\ddot{v}(t)+e^{tJ}J\dot{v}(t)\\
&& {}-\add{x'(t)}{\Bigl(\dot{z}(t)-[e^{tJ}v(t),x'(t)]\Bigr)}\,=\, 0\,.
\end{eqnarray*}
Separating this into $\z$ and $\v$ components yields the desired result.
\end{proof}

\begin{definition}
Let $\g$ denote a geodesic and assume that $\g(t_0)$ is conjugate to 
$\g(0)$ along $\g$. To indicate that the multiplicity of $\g(t_0)$ is $m$,
we shall write $\mcp(t_0)=m$. To distinguish the notions clearly, we shall
denote the multiplicity of $\l$ as an eigenvalue of a specified linear 
transformation by $\mev{\l}$.
\end{definition}

Let $\g$ be a geodesic with $\g(0)=1$ and $\dg(0) = z_0+x_0 \in
\z\ds\v$, respectively, and let $J = J_{z_0}$. If $\g$ is not 
null, we may assume $\g$ is normalized so that $\langle\dg,\dg\rangle = 
\pm 1$. As usual, $\bZ^*$ denotes the set of all integers with 0 removed.
We adapt the usual notation from number theory and write $a|b$ to denote 
that $b$ is a nonzero integral multiple of $a$ for real $a,b$.
\begin{theorem}
Under\label{th1} the preceding assumptions:
\begin{enumerate}
\item if $J=0$ and $x_0=0$, then there are no conjugate points along $\g$;

\item if $J=0$ and $x_0\ne 0$, then $\g(t)$ is conjugate to $\g(0)$
along $\g$ if and only if $-12/t^2$ is an eigenvalue of the linear
operator
$$A : \z\to\z : z\mapsto\left[ x_0,J_z x_0 \right] $$
and $\mcp(t) = \mev(-12/t^2) \le \dim\z$.
\item if $J\ne 0$ and $x_0 = 0$, then $\g(t)$ is conjugate to $\g(0)$
along $\g$ if and only if
$$t \in \bigcup^q_{k=1} \frac{2\pi}{\l_k} \bZ^*,$$
where the $-\l^2_k$ are the negative eigenvalues of $J^2$, and 
$$ \mcp\!\left(\frac{2\pi}{\l_k}n\right) = 
\sum_{\frac{\l_k}{n}\big|\l_h}\!\mev\!\left(-\l_h^2\right) \le \dim\v\, 
.$$

\end{enumerate}
\end{theorem}
If there are no negative eigenvalues in Part 3, then there are of course 
no conjugate points.

Consider the conjugate points from Part 2 of the preceding theorem for
which $z_0=0$.  Denoting this set by $\cal{Z}$, we have $\cal{Z} \subseteq
\exp\v$ (the Lie group exponential map).  Using standard results 
from analytic function theory, we immediately obtain the next result.

\begin{corollary}
This\label{avz} part $\cal{Z}$ of the conjugate locus is an analytic 
variety in the submanifold $\exp\v$ in $N$.\eop
\end{corollary}
For the convenience of the reader, we sketch one way to prove this using
results in the well-known text of Whitney \cite{Wh}.  Complexify.  By
Theorem V.4A, the mapping of coefficients of the characteristic
polynomial of a matrix, thus of the entries, to the eigenvalues is
continuous and proper; {\em cf.}\ the map $\tau$ on p.\,355 in item (D),
second paragraph.  We take $\v$ as our domain, {\em via\/} the composition
of $\tau$ with the map $\{\mbox{components of } x_0\} \to
\{\mbox{coefficients of the characteristic polynomial of } A\}$.  So we
think of the set of eigenvalues as the image of a multifunction of
$x_0\in\v$.  Now the results in Section 1.7 (pp.\,21ff) show that $\tau$
is a holomorphic multifunction, whence so is our composition just defined.
Theorems 4.6F, 5.4A, and/or 7.11B yield analyticity of the image of (a
subset of) the eigenvalues, by the composition of our multifunction with
the mapping defined by the geodesic parametrization.  (Since $\tau$ is
proper, this last composition is at least semiproper, which suffices for
the use of Theorem 7.11B.)  Finally, our conjugate locus is the real part,
thus (real) analytic itself. In the Riemannian case, stronger results 
follow from Warner \cite{Wa}.

\begin{theorem}
In\label{th2} addition, assume that $\dim\z=1$. Then:
\begin{enumerate}
\item if $z_0\ne 0$ and $x_0 = 0$, then $\g(t)$ is conjugate to
$\g(0)$ along $\g$ if and only if
$$t \in \bigcup^q_{k=1} \frac{2\pi}{\l_k}\bZ^* $$
where the $-\l^2_k$ are the negative eigenvalues of $J^2$ and 
$$ \mcp\!\left(\frac{2\pi}{\l_k}n\right) =
\sum_{\frac{\l_k}{n}\big|\l_h}\!\mev\!\left(-\l_h^2\right) \le \dim\v\, 
;$$

\item if $z_0 = 0$ and $x_0\ne 0$, then $\g(t)$ is conjugate to $\g(0)$
along $\g$ if and only if
$$-\frac{12}{t^2} = \ve \langle J_z x_0, J_z x_0 \rangle ,$$
where $z\in\z$ is a unit vector with $\ve = \langle z,z \rangle$, and 
$\mcp(t)=1$;

\item if $z_0\ne 0$ and $x_0\ne 0$, then $\g(t)$ is conjugate to $\g(0)$
along $\g$ if and only if either
$$t \in \bigcup^q_{k=1} \frac{2\pi }{\l_k} \bZ^*$$
where the $-\l^2_k$ are the negative eigenvalues of $J^2$, or $t$ is a
solution of
$$ t\langle Jx_0, \left(e^{-tJ}-I\right)^{-1}x_0\rangle = 
\langle\dg,\dg\rangle $$
in which case $\mcp(t)=1$. For $t=2\pi n/\l_k$, the multiplicities are as 
follows. If $x_0\notin\im(e^{-tJ}-I)$, then 
$$ \mcp(t) = \sum_{\frac{\l_k}{n}\big|\l_h}\!\mev\!\left(-\l_h^2\right) - 
   1\,. $$
If $x_0\in\im(e^{-tJ}-I)$, then choose $v$ such that $(e^{-tJ}-I)v = 
tx_0$ and 
$$ \mcp(t) = \left\{ \begin{array}{cl}
\dps\sum_{\frac{\l_k}{n}\big|\l_h}\!\mev\!\left(-\l_h^2\right) + 1 & 
   \dps\mbox{ if }\, \langle Jx_0,v\rangle = \langle\dg,\dg\rangle\, 
   ,\\[4.5ex]
\dps\sum_{\frac{\l_k}{n}\big|\l_h}\!\mev\!\left(-\l_h^2\right) & 
   \dps\mbox{ if }\, \langle Jx_0,v\rangle \neq \langle\dg,\dg\rangle\,.
   \end{array} \right. $$

\end{enumerate}
\end{theorem}
As we now assume $\z$ is 1-dimensional, Corollary \ref{avz} 
simplifies. (Recall that we are using the Lie group exponential map.)
\begin{corollary}
The\label{smz} part $\cal{Z}$ of the conjugate locus is an analytic 
submanifold of the hypersurface $\exp\v$ in $N$.\eop
\end{corollary}
Alternatively, one could prove this directly from the equation in Part 2.

In the next result, we use the orthogonal direct sum decomposition
$\bigoplus_{j=1}^m \w_j$ where $J$ leaves each $\w_j$ invariant and
each $\w_j$ is an eigenspace of $J^2$. Further, $J^2$ is negative definite
for $1\le j\le q$ and positive definite for $q+1\le j\le m$.
\begin{corollary}
If\label{cor} $J^2$ is moreover diagonalizable, then Parts 2 and 3 can be 
made completely explicit.
\begin{itemize}
\item[2.] if $z_0 = 0$ and $x_0\ne 0$, then $\g(t)$ is conjugate to 
$\g(0)$
along $\g$ if and only if
$$\frac{12}{t^2} = \sum_{l=1}^{m-q} \l^2_{q+l} \ve \langle
B_{q+l},B_{q+l} \rangle - \sum_{k=1}^q \l^2_k \ve \langle A_k,A_k
\rangle ,$$
where $x_0 = \sum_k A_k + \sum_l B_{q+l}$, $A_k \in \w_k$, $B_{q+l} \in
\w_{q+l}$, $-\l^2_k$ and $\l^2_{q+l}$ are the eigenvalues of $J^2$
with $J = J_z$ for $z\in\z$ a unit vector, and $\ve = \langle
z,z \rangle$;

\item[3.] if $z_0\ne 0$ and $x_0\ne 0$, then $\g(t)$ is conjugate to 
$\g(0)$
along $\g$ if and only if either
$$t \in \bigcup^q_{k=1} \frac{2\pi }{\l_k} \bZ^*$$
or $t$ is a solution of
$$ \sum_{k=1}^q \langle A_k,A_k \rangle \frac{\l_k
t}{2}\cot\frac{\l_k t}{2} + \sum_{l=1}^{m-q} \langle B_{q+l},B_{q+l}
\rangle \frac{\l_{q+l} t}{2} \coth\frac{\l_{q+l} t}{2} = 
\langle\dg,\dg\rangle $$
where the $A_h$, $B_i$ are as in the preceding part and $\pm\l_j^2$ are 
the eigenvalues of $J^2$ with $J = J_{z_0}$ now.

\end{itemize}
\end{corollary}

We finish this section by investigating more closely the smoothness of 
this conjugate locus near ({\em i.\,e.,} in a tubular neighborhood of) the
set $\cal{Z}$.

Continue with $\dim\z=1$, $z\in\z$ a unit vector with $\langle
z,z\rangle=\ve$, and $J^2_z$ diagonalizable.  Consider $x_0\in\v$
decomposed as $x_0 = \sum_{k=1}^p A_k + \sum_{l=p+1}^m B_l$ where $A_k$
and $B_l$ are eigenvectors of $J_z^2$ with corresponding eigenvalues
$-\l_k^2,\l_l^2$ where $\l_j>0$ and $\langle x_0,x_0\rangle = \sum_k
\langle A_k,A_k\rangle + \sum_l \langle B_l,B_l\rangle = \pm1$ or 0
according as $x_0$ is unit timelike, unit spacelike, or null,
respectively. Note that the $A_k$, $B_l$, and eigenvalues depend 
analytically on $x_0$.

We shall study 1-parameter families of geodesics whose initial velocities
approach $x_0$.  So for the parameter $a$ consider the geodesic $\g_a$
emanating from $1\in N$ with $\dg_a(0)$ the appropriate linear combination
of $z$ and $x_0$ so that $\g_a$ is unit speed if $x_0$ is nonnull.

$$ \renewcommand{\arraystretch}{1.5}
\begin{array}{c|c|c|c} \dps
\langle x_0,x_0\rangle & \dps\dg_a(0) & a\mbox{ range} & 
   \dps\langle\dg_a,\dg_a\rangle \\\hline
\dps\ve  & \dps\; az+\sqrt{1-a^2}\,x_0\: & \dps\:-1\le a\le 1\: & \dps \ve 
   \\\hline
\dps-\ve & \dps az+\sqrt{1+a^2}\,x_0 & \dps a\in\R & \dps -\ve \\\hline
\dps 0   & \dps az+\sqrt{1+a^2}\,x_0 & \dps a\in\R & \dps a^2\ve        
\end{array} $$

\bigskip
By Corollary \ref{cor}, $\g_a(t)$ is conjugate to $\g_a(0) = 1$ along
$\g_a$ if $t$ solves, for example,
\begin{eqnarray}
(1\mp a^2)\sum_k \langle A_k,A_k \rangle \frac{\l_k
   at}{2}\cot\frac{\l_k at}{2} && \nonumber\\
{} + (1\mp a^2)\sum_l \langle B_{l},B_{l}\rangle \frac{\l_{l} at}{2}
   \coth\frac{\l_{l} at}{2} &=&  \langle \dg_a,\dg_a\rangle \label{1}
\end{eqnarray}
for $a\ne 0$, or if $t$ solves
$$ \frac{12}{t^2} = \sum_{l} \l^2_{l} \ve \langle
   B_{l},B_{l} \rangle - \sum_{k} \l^2_k \ve \langle A_k,A_k \rangle = 
   \Delta^2 $$
for $a=0$.  Since there are no real solutions (hence no conjugate points)
unless the quantity is positive, we have written it as $\Delta^2$ with
$\Delta>0$. Note that the left-hand side of (\ref{1}) and $\Delta$ both 
depend analytically on $x_0$.

We claim that for sufficiently small $a$, equation (\ref{1}) has a 
unique (up to sign), bounded solution $t(a,x_0)$ such that
\begin{equation}
\label{l}
\lim_{a\to 0} t(a,x_0) = 
\pm\frac{2\sqrt{3}}{\Delta } 
\end{equation}
analytically in $x_0$.  We shall only write out the case where $|z| =
|x_0| = 1$ and $t>0$; the other cases can be handled similarly.

First, note that the functions $x\cot x$ and $x\coth x$ have removable 
singularities at 0. Thus we may consider their Maclaurin expansions. To 
simplify notation, we shall supress explicit $x_0$ dependence for a little
while.
\begin{eqnarray*}
\frac{\l_k at}{2}\cot\frac{\l_k at}{2} &=& 1 -
\frac{\left(\l_k at\right)^2}{12} + \, a^4 f_k(a,t)\\
\frac{\l_l at}{2}\coth\frac{\l_l at}{2} &=& 1 +
\frac{\left(\l_l at\right)^2}{12} + \, a^4 f_l(a,t)
\end{eqnarray*}
Here $a^4 f_n(a,t)$ denote the remainders obtained by summing the rest of
the series, with $f_n(0,t) \neq 0$.  Next, substitute these into equation
(\ref{1}) with $\langle \dg_a,\dg_a\rangle = 1$.
\begin{eqnarray*}
f(a,t) &=&(1-a^2)\sum_k \langle A_k,A_k \rangle \left[1 - \frac{\left(\l_k
   at\right)^2}{12} + \, a^4 f_k(a,t) \right] \nonumber\\
&& {} + (1-a^2)\sum_l \langle B_l,B_l \rangle \left[ 1 + \frac{\left(\l_l
   at\right)^2}{12} + \, a^4 f_l(a,t) \right] - 1 \,=\, 0
\end{eqnarray*}
Using $\sum_k \langle A_k,A_k\rangle + \sum_l \langle B_l,B_l\rangle =1$, 
we obtain
\begin{equation}
\label{s}
a^2\left( 1 - \frac{t^2}{12}\Delta^2 + a^2 R(a,t) \right) = a^2 g(a,t)
   = 0\, ,
\end{equation}
where $a^2 R(a,t)$ is the combination of all remaining terms with $R(0,t) 
\neq 0$.

Observe that $g(0,2\sqrt{3}/\Delta) = 0$ and 
$$ \frac{\partial g}{\partial t}\left( 0,\frac{2\sqrt{3}}{\Delta}\right) 
= -\frac{\Delta}{\sqrt{3}} \neq 0\,.$$
Thus we may apply the Implicit Function Theorem to conclude that, near 
$(0,2\sqrt{3}/\Delta)$, there exists a unique, positive, analytic solution 
$t=t(a,x_0)$. It follows from (\ref{s}) that this $t$ satisfies (\ref{l}). 
Now we use Corollary \ref{smz} to conclude the proof of the next result.
(Recall that we are using the Lie group exponential map.)
\begin{proposition}
The\label{scl} conjugate locus is smooth (indeed, analytic) across (or
through) the hypersurface $\exp\v$ in $N$.
\end{proposition}
\begin{proof}
Let $n=\dim N$ and consider $p\in\cal{Z}$. Since $\cal{Z}$ is smooth, 
there are $n-2$ smooth coordinate functions centered at $p$ comprising a 
chart in $\cal{Z}$. Now $t=t(a,x_0)$ provides another independent, smooth 
coordinate function centered at $p$. Together with the $n-2$ we already 
have, we now have a chart in the conjugate locus centered at $p$. Thus 
near $\cal{Z}$, the conjugate locus is a smooth, codimension-1 
submanifold.
\end{proof}
Note that this applies only in a (possibly very small) tubular
neighborhood of $\cal{Z}$ in $N$.  Away from $\cal{Z}$ (outside this
tubular neighborhood), the conjugate locus is merely an analytic variety
in $N$.  In the Riemannian case, stronger results follow from Warner
\cite{Wa}.

\section{Proof of Theorem \protect\ref{th1}}
\label{pf1}

The first part is an immediate consequence of Proposition \ref{4.18} (4.18
in \cite{CP4}).

For the second part, let $Y(t) = z(t) + e^{tJ}v(t)= z(t) + v(t)$ be a
nontrivial Jacobi field along $\g$ with $Y(0)=Y(t_0)=0$ for some $t_0 \neq
0$.  Then by Proposition \ref{yj} and $J=0$, we have
\begin{eqnarray}
\dot{z}(t) -[ v(t) , x_0 ] &=& \zeta\, ,\label{jp1}\\
\ddv(t) - J_{\zeta}x_0 &=& 0\, .\label{jp2}
\end{eqnarray}
for some constant vector $\zeta\in\z$. By (\ref{jp2}) and $ v(0)=v(t_0)=0$,
\begin{equation}
\label{jp3}
v(t) = \half t(t-t_0)J_{\zeta}x_0\, .
\end{equation}
Substituting this into (\ref{jp1}), integrating, and using $z(0)=0$ yields
\begin{equation}
\label{jp4}
z(t) = (\sixth t^3 - \quar t_0 t^2 )[J_{\zeta}x_0 , x_0] +t\zeta\,.
\end{equation}
Using $z(t_0)=0$ on (\ref{jp4}) we find
\begin{equation}
\label{jp5}
[x_0 , J_{\zeta}x_0 ] = -\frac{12}{t_0^2}\zeta\,.
\end{equation}
If $\zeta=0$, then (\ref{jp3}) and (\ref{jp4}) imply that $Y(t) = 0$ for 
all $t$, contradicting the nontriviality of $Y$.  Thus $\zeta\neq 0$ whence
(\ref{jp5}) implies that $-12/t_0^2$ is an eigenvalue of the linear map
$A:\z \to \z : z\mapsto [x_0 , J_z x_0]$.  Conversely, if $-12/t_0^2$
is an eigenvalue of $A$ there is an eigenvector $\zeta \in \z$ satisfying
(\ref{jp5}).  Then $Y(t) = z(t) + v(t)$, where $z(t)$ and $v(t)$ are given
by (\ref{jp4}) and (\ref{jp3}), respectively, is a nontrivial Jacobi field
along $\g$ with $Y(0)=Y(t_0)=0$. The multiplicity is obvious.

\medskip

We turn now to the proof of the third part of Theorem \ref{th1}. Assume
that $\g(t_0)$ is conjugate to $\g(0)$ along $\g$. Then there exists a
nonzero Jacobi field $Y(t) = z(t) + e^{tJ}v(t)$ along $\g$ with
$Y(0)=Y(t_0)=0$.  By Proposition \ref{yj} and $x_0 =0$ we have
\begin{eqnarray}
\ddot{z}(t) &=& 0\, ,\label{jp6}\\
\ddv + J\dot{v}&=&0\,.\label{jp7}
\end{eqnarray}
{F}rom (\ref{jp6}) and $z(0)=z(t_0)=0$, we see that $z(t) = 0$. 

Let
$$p(t)=\prod_{j=1}^p \left( t^2-2a_j t+\left( a_j^2+b_j^2 \right)
   \right)^{k_j} \prod_{j=1}^q \left( t^2 +\lambda_j^2 \right)^{l_j}
   \prod_{j=1}^r \left( t-\mu_j \right)^{m_j}\cdot t^s $$
be the characteristic polynomial of $J$ where $a_j\neq 0$, $b_j\neq 0$,
$\lambda_j > 0$, and $\mu_j\neq 0$ are all real.  Then equation
(\ref{jp7}) can be split into
\begin{eqnarray}
\ddot{v}_1 + J \dot{v}_1 &=& 0\, , \label{jp8}\\
\ddot{v}_2 + J\dot{v}_2 &=&0\, , \label{jp9}
\end{eqnarray}
where $ v = v_1 + v_2$ with $\v_1 \in \k$ for
$$ \k = \ker\prod_{j=1}^p \left( J^2-2a_j J+ \left( a_j^2+b_j^2
   \right)I \right)^{k_j} \prod_{j=1}^q \left( J^2 +\lambda_j^2I
   \right)^{l_j} \prod_{j=1}^r \left( J-\mu_jI \right)^{m_j} $$
and $ v_2 \in\ker J^s$.  Applying $J^{s-1}$ to both sides of (\ref{jp9}),
we have $ J^{s-1}\ddot{v}_2 = 0$.  This and $J^{s-1}v_2(0) =
J^{s-1}v_2(t_0) = 0$ imply that $J^{s-1}v_2 = 0$.  Next, applying
$J^{s-2}$ we have $J^{s-2}\ddot{v}_2 = 0$, which similarly implies that
$J^{s-2}v_2 = 0$.  Continuing this process, we conclude that $v_2 = 0$. 
Note that $J$ is invertible on $\k$, essentially because $\ker 
J\subseteq\ker J^s$; see \cite{JP2} for some more detailed computations of 
this type.
\begin{lemma}
Let\label{ex} $A$ be a real nonsingular matrix and $0\neq t\in\R$.
Then $\ker\left(e^{tA}-I\right) = \bigoplus_{t\l_j \in 2\pi\bZ^*}
\ker(A^2+\l_j^2 I)$ where the $-\l_j^2$ are the negative eigenvalues (if 
any) of $A^2$.
\end{lemma}
\begin{proof}
Let $v \in \ker(A^2 + \l_j^2 I)$ with $t\l_j \in 2\pi\bZ^*$. Then $e^{tA}v
= (\cos t\l_j)v + (\sin t\l_j)Jv/\l_j = v$ so $\left(e^{tA}-I\right)v = 0$
and $v \in \ker\left(e^{tA}-I\right)$. 

The converse follows {\em via\/} a straightforward computation using the 
Jordan canonical form of $A$; {\em cf.} \cite{JP2} for similar work.
\end{proof}
Using $v_1(0)=0$, from (\ref{jp8}) we get $v_1(t) = (e^{-tJ}-I)v$ for some
constant vector $0\neq v\in\k$.  Since $v_1(t_0)=0$, we conclude from this
and the lemma that $t_0 \in \bigcup_{j=1}^q (2\pi/\l_j) \bZ^*$ and
$$ v \in \ker\left(e^{t_0J}-I\right)\cap\k = \bigoplus_{t_0\l_j \in 
   2\pi\bZ^*}\!\ker(J^2 + \l_j^2 I)\,. $$

Moreover, any Jacobi field vanishing at 0 and $t_0$ is easily seen to be 
of the form
$$ Y(t) = \left( e^{-tJ}-I\right)v \quad\mbox{ for some }\quad v \in 
   \bigoplus_{t_0\l_j \in 2\pi\bZ^*} \!\ker(J^2 + \l_j^2 I)\, ,$$
so the mulitplicity follows by noting that $t_0\l_h \in
2\pi\bZ^*$ for $t_0 = 2\pi n/\l_k$ is equivalent to $\l_h = m\l_k/n$ for
some nonzero integer $m$.

\section{Proof of Theorem \protect\ref{th2}}
\label{pf2}

Again the first part is immediate, this time by the third part of Theorem
\ref{th1}. The second part follows from the second part of Theorem 
\ref{th1} upon noting that $A$ is now a $1\times 1$ matrix and that 
$$ \langle [x_0,J_z x_0],z \rangle = \langle J_z x_0,J_z x_0 
\rangle .$$

For the third part, let $Y(t) = \alpha (t) z_0 + e^{tJ}v(t)$ be a Jacobi
field along $\g$ with $Y(0) = Y(t_0) = 0$.  Then by Proposition \ref{yj}
we have
\begin{eqnarray*}
\da(t)z_0 -[e^{tJ}v(t) , e^{tJ}x_0 ] &=& cz_0\, ,\\
e^{tJ}\ddv(t) + e^{tJ}J\dv(t) -cJe^{tJ}x_0 &=& 0
\end{eqnarray*}
for a constant $c$.  Note that $[e^{tJ}x, e^{tJ}y]=[x,y]$ since $\dim \z =
1$, and recall that $J$ commutes with $e^{tJ}$.  Using these, we simplify
the preceding equations as
\begin{eqnarray}
\da(t)z_0 -[v(t) , x_0 ]& =& cz_0\, ,\label{jp10} \\
\ddv + J\dot{v} -cJx_0 &=& 0\,.\label{jp11}
\end{eqnarray}
Taking the inner product of (\ref{jp10}) with $z_0$ and rearranging,
\begin{equation}
\label{jp12}
\da +\frac{\langle Jx_0 , v \rangle }{\langle z_0 , z_0 \rangle } = c
\end{equation}
and the general solution of (\ref{jp11}) satisfying $v(0)=0$ is 
\begin{equation}
\label{gs}
v(t) = \left( e^{-tJ}-I \right) v_0 + ctx_0
\end{equation}
for some constant vector $v_0\in\v$.

\bigskip
Assume that $\g(t_0)$ is conjugate to $\g(0)$ and
\begin{equation}
\label{hyp1}
t_0 \in \bigcup_{k=1}^q \frac{2\pi }{\l_k}\bZ^* .
\end{equation}
We distinguish two cases according as $x_0$ is or is not in the image of 
$e^{-t_0 J}-I$.

\medskip
If $x_0 \notin \im\left( e^{-t_0 J}-I \right)$, then $v(t_0)=0$ and 
(\ref{gs}) imply $v_0 \in \ker\left( e^{-t_0 J}-I \right)$ and $c=0$. 
Using (\ref{jp12}) and integrating under $\alpha(0)=0$,
$$ \alpha(t) + \frac{\langle Jx_0, (-J)^{-1}\left( e^{-tJ}-I 
   \right)v_0-tv_0 \rangle }{\langle z_0,z_0\rangle } = 0\,. $$
Since $\alpha(t_0)=0$, it follows that $\langle Jx_0,v_0\rangle = 0$. Thus
if $Y$ is a Jacobi field along $\g$ such that $Y(0)=Y(t_0)=0$, then 
(\ref{hyp1}) implies that
$$ v(t) = \left( e^{-tJ}-I \right)v_0\, ,$$
for some $v_0 \in \ker\left( e^{-t_0 J}-I \right)$ with $\langle 
Jx_0,v_0\rangle = 0$, and
$$ \alpha(t) = \frac{\langle Jx_0,tv_0 - (-J)^{-1}\left( e^{-tJ}-I \right)
v_0\rangle }{\langle z_0,z_0\rangle }\,. $$
\begin{lemma}
For\label{limg} $x\in\v$, $\langle x,v\rangle = 0$ for all $v \in
\ker\left( e^{-t_0 J}-I \right)$ if and only if $x \in \im\left( e^{-t_0
J}-I \right)$.
\end{lemma}
We omit the proof, which is a reasonably straightforward computation in
linear algebra.  One does need to note that $x_0 \notin \im\left( e^{-t_0
J}-I \right)$ if and only if $Jx_0 \notin \im\left( e^{-t_0 J}-I \right)$,
to use the orthogonal decomposition $\v = \left( \bigoplus_{t_0\l_i \in 
2\pi\bZ} \ker(J^2+\l_i^2)^{l_i}\right) \ds \v'$ where $\v'$ is a 
$J$-invariant, nondegenerate subspace of $\v$ on which $\left( e^{-t_0 
J}-I \right)$ is invertible, and to use an appropriate block-diagonal 
representation of $J$ on subspaces of the kernel summands. See, for 
example, \cite{JP2} for several similar such computations.

This lemma now yields $\mcp(t_0) = \dim\ker\left( e^{-t_0 J}-I \right) - 
1$ and the multiplicity formula follows as previously.

\medskip
If $x_0 \in \im\left( e^{-t_0 J}-I \right)$, then we can take $v_1 \in \v$
such that $\left( e^{-t_0 J}-I \right)v_1 = t_0x_0$.  Using (\ref{gs}) and
$v(t_0)=0$, we obtain $v_0 = v_2-cv_1$ for some $v_2 \in \ker\left(
e^{-t_0 J}-I \right)$.  From (\ref{jp12}), integrating under $\alpha(0)=0$
yields
$$ \alpha(t) + \frac{\langle Jx_0, (-J)^{-1}\left( e^{-tJ}-I 
   \right)v_0-tv_0\rangle }{\langle z_0,z_0\rangle } = ct\,. $$
Using $\alpha(t_0)=0$, Lemma \ref{limg}, and the skewadjointness of $J$, a 
short computation shows
\begin{equation}
\label{nc}
c\left(\langle Jx_0,v_1\rangle - \langle\dg,\dg\rangle\right)=0\,.
\end{equation}
By Lemma \ref{limg}, $\langle Jx_0,v_1\rangle$ is independent of the
choice of such a $v_1$.  We distinguish two subcases according as $\langle
Jx_0,v_1\rangle$ is or is not equal to $\langle\dg,\dg\rangle$.

\smallskip
If $\langle Jx_0,v_1\rangle = \langle\dg,\dg\rangle$, then we have
$$ v(t) = \left( e^{-tJ}-I \right)(v_2-cv_1) + ct_0x_0\, ,$$
$$ \alpha(t) = ct - \frac{\langle Jx_0, (-J)^{-1}\left( e^{-tJ}-I \right) 
   (v_2-cv_1) + tcv_1\rangle }{\langle z_0,z_0\rangle } $$
for $v_2 \in \ker\left( e^{-t_0 J}-I \right)$ and arbitrary scalar $c$, 
whence 
$$\mcp(t_0) = \dim\ker\left( e^{-t_0 J}-I \right) + 1$$
and the multiplicity formula follows as before.

\smallskip
If $\langle Jx_0,v_1\rangle \neq \langle\dg,\dg\rangle$, then (\ref{nc}) 
implies $c=0$ and we find
$$ v(t) = \left( e^{-tJ}-I \right)v_0 $$
for some $v_0 \in \ker\left( e^{-t_0 J}-I \right)$,
$$ \alpha(t) = -\frac{\langle x_0,\left( e^{-tJ}-I \right)v_0\rangle 
   }{\langle z_0,z_0\rangle }\,, $$
and $\mcp(t_0) = \dim\ker\left( e^{-t_0 J}-I \right)$, from which the 
desired formula follows as before.

\bigskip
Now assume that $\g(t_0)$ is conjugate to $\g(0)$ and
\begin{equation}
\label{hyp}
t_0 \notin \bigcup_{k=1}^q \frac{2\pi }{\l_k}\bZ^* .
\end{equation}
If $c=0$, then $v(t_0)=0$ implies $\left( e^{-t_0 J}-I \right) v_0
= 0$ in (\ref{gs}).  But then (\ref{hyp}) and Lemma \ref{ex} forces
$v_0=0$ in (\ref{gs}).  Thus $u=0$ and $\langle Jx_0,u\rangle = 0$, so
together with (\ref{jp12}) and $\alpha(0) = \alpha(t_0) = 0$ these imply
that $\alpha = 0$.  Consequently, $Y=0$ which contradicts our assumption
that $Y$ is a nontrivial Jacobi field.  Therefore $c\neq 0$.

Using (\ref{gs}) and $v(t_0) = 0$, we get
\begin{equation}
\label{3j1.5}
v(t) = ctx_0 - \left( e^{-tJ}-I \right) \left(e^{-t_0 J}-I 
\right)^{-1}ct_0 x_0 \,.
\end{equation}
Substituting this into (\ref{jp12}), we find
$$ \da (t) - \frac{\langle Jx_0, \left( e^{-tJ}-I \right) \left(e^{-t_0 
   J}-I\right)^{-1}ct_0 x_0 \rangle }{\langle z_0,z_0 \rangle } = c\,. $$
Integrating with $\alpha (0)=0$ yields
\begin{eqnarray}
\lefteqn{\alpha (t) =  ct }\label{3j1.6}\\
&&{}+\frac{\langle Jx_0, (-J)^{-1} \left( e^{-tJ}-I 
   \right) \left(e^{-t_0 J}-I\right)^{-1}ct_0 x_0 - t \left( e^{-t_0 J}-I 
   \right)^{-1}ct_0 x_0 \rangle }{\langle z_0,z_0 \rangle }\,.\nonumber
\end{eqnarray}
Since $\alpha(t_0)=0$, from (\ref{3j1.6}) we obtain
$$ 0 = 1 + \frac{ \langle Jx_0, (-J)^{-1}x_0 -t_0\left( e^{-t_0 J}-I 
   \right)^{-1}x_0 \rangle }{\langle z_0,z_0 \rangle }\,.$$
Thus
$$ 0 = -\langle z_0,z_0\rangle - \langle x_0,x_0\rangle + t_0\langle Jx_0,
   \left( e^{-t_0 J}-I\right)^{-1}x_0 \rangle $$
or
\begin{eqnarray*}
\langle\dg,\dg\rangle &=& \langle z_0,z_0\rangle + \langle x_0,x_0\rangle \\
&=& t_0 \langle Jx_0, \left( e^{-t_0 J}-I \right)^{-1}x_0 \rangle\,.
\end{eqnarray*}
Therefore, if $\g(t_0)$ is conjugate to $\g(0)$ along $\g$ and (\ref{hyp})
holds, then $t_0$ is a solution of
\begin{equation}
\label{3j1.7}
 t\langle Jx_0, \left( e^{-tJ}-I\right)^{-1}x_0 \rangle = 
   \langle\dg,\dg\rangle\,.
\end{equation}
The multiplicity follows from the fact that $v$ in (\ref{3j1.5}) and 
$\alpha$ in (\ref{3j1.6}) are uniquely determined by the scalar $c$.

\section{Proof of Corollary \protect\ref{cor}}
\label{pf3}

For Part 2, let $z$ be a unit vector in $\z$ and consider the inner
product $\langle [ x_0 , J_z x_0],z \rangle$.  In this part of the
proof $J$ denotes $J_{z}$.  Thus,
\begin{eqnarray*}
\langle [x_0,Jx_0],z \rangle &=& \langle Jx_0,Jx_0 \rangle\\
   &=& \langle -J^2x_0,x_0 \rangle .
\end{eqnarray*}
Since $J^2$ is nonsingular and diagonalizable, $\v = \bds_{k=1}^q \w_k \ds
\bds_{l=1}^{m-q} \w_{q+l}$ where $\w_k$ is the eigenspace of $-\l_k^2$ and
$\w_{q+l}$ is the eigenspace of $\l_{q+l}^2$.  Writing $x_0 = \sum_{k=1}^q
A_k + \sum_{l=1}^{m-q} B_{q+l}$ with $A_k \in \w_k$ and $B_{q+l} \in
\w_{q+l}$, we have
\begin{eqnarray*}
\langle [x_0,Jx_0],z \rangle &=& \langle -J^2x_0,x_0 \rangle \\
   &=& \sum_{k=1}^q \l_k^2 \langle A_k,A_k \rangle - \sum_{l=1}^{m-q}
       \l_{q+l}^2 \langle B_{q+l},B_{q+l} \rangle .
\end{eqnarray*}
This implies that
$$ [x_0,Jx_0] = \ve\left( \sum_{k=1}^q \l_k^2 \langle A_k,A_k \rangle - 
   \sum_{l=1}^{m-q} \l_{q+l}^2 \langle B_{q+l},B_{q+l} \rangle \right) z
$$
where $\ve = \langle z,z \rangle$. By the second part of Theorem
\ref{th2}, $\g(t_0)$ is conjugate to $\g(0)$ along $\g$ if and only if
$$ -\frac{12}{t_0^2} = \ve\left( \sum_{k=1}^q \l_k^2 \langle A_k,A_k
   \rangle - \sum_{l=1}^{m-q} \l_{q+l}^2 \langle B_{q+l},B_{q+l} \rangle
   \right) .$$

\medskip

Finally, we give the proof of Part 3 of Corollary \ref{cor}.  Now we have
$\dg(t) = z_0 + e^{tJ}x_0$ with $z_0\neq 0 \neq x_0$ and $J =
J_{z_0}$ again, as usual.  We shall continue with the notations
for eigenvalues, eigenspaces, and components of $x_0$ as in the preceding
Part 2.

It will suffice to show that
\begin{eqnarray*}
\lefteqn{t\langle Jx_0, \left(e^{-tJ}-I\right)^{-1}x_0\rangle = }\\
&& \sum_{k=1}^q \langle A_k,A_k \rangle \frac{\l_k
   t}{2}\cot\frac{\l_k t}{2} + \sum_{l=1}^{m-q} \langle B_{q+l},B_{q+l}
   \rangle \frac{\l_{q+l} t}{2} \coth\frac{\l_{q+l} t}{2}\, .
\end{eqnarray*}
Using the $J$-invariance of each $\w_j$ and the orthogonality of the 
direct sum decomposition, we have
\begin{eqnarray*}
\lefteqn{t\langle Jx_0, \left(e^{-tJ}-I\right)^{-1}x_0\rangle = }\\
&& t\sum_{k=1}^q\langle JA_k, (e^{-tJ}-I)^{-1}A_k\rangle + 
t\sum_{l=1}^{m-q} \langle JB_{q+l}, (e^{-tJ}-I)^{-1}B_{q+l}\rangle .
\end{eqnarray*}
If $A_k\neq 0$, then $\lsp A_k, \l_k^{-1}JA_k\rsp$ is a 2-dimensional
complement in $\w_k$ of a $J$-invariant subspace $\w'_k$. On this 
planar subspace, the matrix of $e^{-tJ}-I$ is
$$ \left[ \begin{array}{cc} \cos t\l_k -1 & \sin t\l_k \\
   -\sin t\l_k & \cos t\l_k -1 \end{array} \right] $$
so the matrix of $(e^{-tJ}-I)^{-1}$ is
$$ \left[ \begin{array}{cc} -\half\csc^2(t\l_k/2) & -\half\cot(t\l_k/2) 
   \raisebox{0ex}[2.7ex][0ex]{} \\[1ex]
   \half\cot(t\l_k/2) & -\half\csc^2(t\l_k/2) \raisebox{0ex}[0ex][1.7ex]{}
   \end{array} \right] $$
and it follows that
\begin{eqnarray}
\langle JA_k, \left(e^{-tJ}-I\right)^{-1}A_k\rangle &=& \langle JA_k,
   -\half\csc^2\frac{t\l_k}{2}A_k + \half\cot\frac{t\l_k}{2}
   \frac{JA_k}{\l_k} \rangle \nonumber\\
&=& \frac{\l_k}{2}\cot\frac{t\l_k}{2}\langle A_k,A_k\rangle \label{3j2.1}
\end{eqnarray}

Now assume $JB_{q+l} \neq\pm \l_{q+l}B_{q+l}$. Then $B_{q+l}$ and 
$JB_{q+l}$ are linearly independent, so we may consider the plane 
$\lsp\l_{q+l}^{-1}JB_{q+l}+B_{q+l}, \l_{q+l}^{-1}JB_{q+l}-B_{q+l}\rsp$ in 
$\w_{q+l}$. Here the matrix of $e^{tJ}-I$ is
$$ \left[ \begin{array}{cc} \dps e^{-t\l_{q+l}}-1 & 0 
   \raisebox{0ex}[2.6ex][0ex]{} \\[1ex]
   0 & \dps e^{t\l_{q+l}}-1 \end{array} \right] $$
so the matrix of $(e^{-tJ}-I)^{-1}$ is
$$ \left[ \begin{array}{cc} \dps\frac{1}{e^{-t\l_{q+l}}-1} & 0 
   \raisebox{0ex}[4ex][0ex]{} \\[2ex]
   0\raisebox{0ex}[0ex][2.6ex]{} & \dps\frac{1}{e^{t\l_{q+l}}-1} 
\end{array} \right] $$
and we obtain
\begin{eqnarray*}
\lefteqn{(e^{-tJ}-I)^{-1}B_{q+l} = }\\[1ex]
&& \frac{\half}{e^{-t\l_{q+l}}-1}\left( \frac{1}{\l_{q+l}}JB_{q+l} +
   B_{q+l} \right) - \frac{\half}{e^{t\l_{q+l}}-1} \left(
   \frac{1}{\l_{q+l}} JB_{q+l} - B_{q+l} \right)
\end{eqnarray*}
whence
\begin{eqnarray}
\langle JB_{q+l}, (e^{-tJ}-I)^{-1}B_{q+l} \rangle  &=& 
   \frac{\l_{q+l}}{2}\langle B_{q+l},B_{q+l}\rangle \left(
   \frac{1}{e^{t\l_{q+l}}-1} - \frac{1}{e^{-t\l_{q+l}}-1} \right) 
   \nonumber\\[1ex]
&=& \frac{\l_{q+l}}{2}\langle B_{q+l},B_{q+l}\rangle 
   \coth\frac{t\l_{q+l}}{2}\,.\label{3j2.2}
\end{eqnarray}

Finally, assume $JB_{q+l} = \pm\l_{q+l}B_{q+l}$. Then
\begin{eqnarray*}
(e^{-tJ}-I)B_{q+l} &=& (e^{\mp t\l_{q+l}}-1)B_{q+l}\, ,\\
(e^{-tJ}-I)^{-1}B_{q+l} &=& \frac{1}{e^{\mp t\l_{q+l}}-1}B_{q+l}\, .
\end{eqnarray*}
Using $\langle JB_{q+l},B_{q+l} \rangle =0$, it follows much as before 
that
\begin{equation}
\label{3j2.3}
\langle JB_{q+l}, (e^{-tJ}-I)^{-1}B_{q+l} \rangle = \frac{\l_{q+l}}{2} 
\langle B_{q+l},B_{q+l} \rangle \coth\frac{t\l_{q+l}}{2}\,.
\end{equation}
Combining (\ref{3j2.1}), (\ref{3j2.2}), and (\ref{3j2.3}) we obtain
\begin{eqnarray*}
\lefteqn{t\langle Jx_0, (e^{-tJ}-I)^{-1}x_0 \rangle }\\
&=& t\sum_{k=1}^q \langle JA_k,
   (e^{-tJ}-I)^{-1}A_k \rangle + t\sum_{l=1}^{m-q} \langle JB_{q+l}, 
   (e^{-tJ}-I)^{-1}B_{q+l} \rangle \\
&=& \sum_{k=1}^q \langle A_k,A_k \rangle \frac{t\l_k}{2}
\cot\frac{t\l_k}{2} + \sum_{l=1}^{m-q} \langle B_{q+l},B_{q+l} \rangle 
\frac{t\l_{q+l}}{2} \coth\frac{t\l_{q+l}}{2}
\end{eqnarray*}
as desired.

\end{document}